\documentclass[reqno]{amsart}
\usepackage{amssymb}

\begin{document}
\title{}
%\author{S. Abbas}
\title[\hfilneg \hfil Fractional Optimal Control]
{Approximate controllability and optimal control of impulsive fractional semilinear delay differential equations with non-local conditions}

\author[L.Mahto, S. Abbas \hfil \hfilneg]
{Lakshman Mahto, Syed Abbas${^*}$, \\ School of Basic Sciences
\\ Indian Institute of Technology Mandi \\ Mandi, H.P. - 175001, India.}

\date{}
 \keywords{Fractional differential
equation, fixed point theorem, controllability, optimal control, lower semi-continuity \\ $*:$ Corresponding author,
Email: sabbas.iitk@gmail.com, \ abbas@iitmandi.ac.in}

\begin{abstract} In this paper we study the approximate controllability and existence of optimal control of impulsive fractional semilinear delay differential equations with non-local conditions. We use Sadovskii's fixed point theorem, semigroup theory of linear operators and direct method for minimizing a functional to establish our results. At the end we give an example to illustrate our analytical findings.
\end{abstract}

\maketitle \numberwithin{equation}{section}
\newtheorem{theorem}{Theorem}[section]
\newtheorem{lemma}[theorem]{Lemma}
\newtheorem{example}[theorem]{Example}
\newtheorem{proposition}[theorem]{Proposition}
\newtheorem{corollary}[theorem]{Corollary}
\newtheorem{remark}[theorem]{Remark}
\newtheorem{definition}[theorem]{Definition}
\newtheorem{assumption}[theorem]{Assumption}

   \vskip .5cm \noindent
   \section{ \textbf{Introduction}}

In this work we consider the following impulsive delay differential equations of fractional order $\alpha \in (1,2)$ with non-local conditions,
\begin{eqnarray}
D_t^{\alpha}x(t)&=& Ax(t)+D_t^{\alpha-1}(f(t,x(t),x_t)+Bu(t)), \quad t \in
[0,T], t\neq t_k \nonumber \\ \Delta x(t)|_{t=t_k}&=& I_k(x(t_k^-)),\quad k=1,2,\cdots ,K \nonumber \\
x(s)+g(x_t)(s)&=& \phi(s) \quad s \in [-r,0] \label{meq}
\end{eqnarray}
where $f:[0,T]\times X\times \mathcal{C}_r\rightarrow X, \mathcal{C}_r=PC([-r,0],X), g:\mathcal{C}_r\rightarrow X, I_{k} : X\rightarrow X, B:U \rightarrow X$  is bounded linear operator on separable reflexive Hilbert space $U$ with norm $\|\cdot\|_U$ and the operator $A:D(A)\subset X \rightarrow X$ is a linear densely defined of sectorial type in a separable Banach space $X$ with norm $\|\cdot\|_X.$  Recently fractional differential equations attract many mathematicians and scientists because of their usefulness in the various problems coming from engineering and physical science. These kind of equations are generalization of ordinary differential equations to arbitrary non integer orders. The origin of fractional calculus goes back to Newton and Leibniz in the seventieth century.  It is widely and efficiently used to describe many phenomena arising in engineering, physics, economy and allied sciences. Recent investigations have shown that many physical systems can be represented more accurately through fractional derivative formulation \cite{pod}. Fractional differential equations, therefore find numerous applications in the field of visco-elasticity, feed back amplifiers, electrical circuits, control theory, electro analytical chemistry, fractional multipoles, neuron modelling encompassing different branches of physics, chemistry and biological sciences \cite{pod1}. Many physical processes appear to exhibit fractional order behavior that may vary with time or space. The fractional calculus has allowed the operations of integration and differentiation to any arbitary order. The order may take on any real or imaginary value. The existence and uniqueness of the solutions of fractional differential equations have been shown by many authors
\cite{abbas, chen, kilbas, lak2, lizama, pod1, salem}.

In mathematical control theory controllability is one of the important concept in which one study  the steering of dynamical system from given initial state to any other state or in neighbourhood of the state under some admissible control input. Many authors \cite{chang, chen, lizama, mohmudov1, mohmudov2,  shen, wang1,  xiang} have been studied controllability of the semilinear evolution equations. Nonlocal condition is more realistic as compared to classical condition in physical systems,  where it involves dependent of initial condition on further states as well. There are several  study \cite{abbas, chen} on controllability of semilinear differential equations with nonlocal initial conditions in Banach spaces.

Impulsive differential equation provides a realistic framework of modeling systems in fields like  population dynamics, control theory, physics, biology and medicine, when the dynamics undergo some abrupt changes at certain moments of time like earthquake, harvesting, shock and so forth. Milman and  Myshkis \cite{milman} first introduced impulsive differential equations in 1960. Followed by their work (Milman and Myshkis), there are several monographes and papers written by many authors like Bainov and Simeonov \cite{bainov1}, Benchohra et al. \cite{benchohra}, Lakshmikantham et al. \cite{laksh}, Samoilenko and Perestyuk \cite{samoilenko} and Mahto et al. \cite{mahto}. Controllability of impulsive semilinear differential equations with nonlocal conditions has been studied in \cite{chen}. Controllability of impulsive fractional evolution equations with nonlocal conditions has been studied in \cite{ amar}.

In several fields like biology, population dynamics and so forth problems with hereditary are best modelled by delay differential equations \cite{hale}. So the problems with impulsive effects and hereditary properties could be modelled by impulsive delay differential equations. Controllability of impulsive semilinear differential equations with finite delay have been studied in \cite{amar,selvi}.

Very little work \cite{ozdemir, wang2, xiang} has been done in optimal control fractional order system. Wang et al.  \cite{wang1} studied fractional finite time delay evolution systems and optimal controls in infinite-dimensional spaces.

%Cueves and Lizama in \cite{lizama} have shown the existence of
%almost automorphic solution of the problem (\ref{meq}) for $t\in
%\mathbb{R}.$
Motivated by the above studied, in this work we establish the approximate controllability and existence of optimal control of impulsive fractional semilinear delay differential equations with nonlocal conditions. In section 3, to establish approximate controllability of the system (\ref{meq}), we prove the existence of mild solution of problem (\ref{meq}) using Sadovskii's fixed point theorem. In section 4, we prove the existence of optimal control of the system (\ref{meq}). At the end  we  give an example to illustrate our analytical results. The results proven in this manuscript are new and interesting.

\noindent
 \section{\textbf{Preliminaries}}
Denote $B(X)$ be the Banach space of all linear and bounded operators on $X$ endowed with the norm $\|\cdot\|_{B(X)}$ and $\mathcal{C}=\mathcal{C}(\mathbb{R},X)$ the set of all continuous functions from $\mathbb{R}$ to $X.$

\textbf{Sectorial operator:} A closed and linear operator $A$ is said to be sectorial of type $\omega$ and angle $\theta$ if there exists $0 <\theta < \frac{\pi}{2}, M_1> 0$ and $\omega \in \mathbb{R}$ such that its resolvent exists outside the sector $$\omega+S_{\theta}:= \{\omega+\lambda: \lambda \in \mathbb{C}, \ |arg(-\lambda)| < \theta\},$$ and $$\|(\lambda- A)^{-1}\| \le \frac{M_1}{|\lambda-\omega|}, \ \lambda \not \in \omega+ S_{\theta}.$$ Sectorial operators are well studied in the literature. For a recent reference including several examples and
properties we refer the reader to \cite{haase}. Note that an operator A is sectorial of type $\omega$ if and only if $\lambda I - A$ is sectorial of type $0.$

The equation (\ref{meq}) can be thought of a limiting case \cite{lizama} of the following equation
\begin{eqnarray}
z^{\prime}(t)=\int_0^t
\frac{(t-s)^{\alpha-2}}{\Gamma(\alpha-1)}Az(s)ds+f(t,x(t),x_t), \ t
\ge 0, \quad z(0)=x_0, \label{meq1}
\end{eqnarray}
in the sense that the solutions are asymptotic to each other as $t \rightarrow \infty.$ If we consider the operator $A$ is sectorial of type $\omega$ with $\theta \in [0,\pi(1-\frac{\alpha}{2})),$ then the problem (\ref{meq1}) is well posed \cite{custa}. Thus we can use variation of parameter formulae to get
$$z(t)=S_{\alpha}(t)(\phi(0)-g(x)(0))+\int_{0}^tS_{\alpha}(t-s)f(s,x(s),x_s)ds, \quad t \ge 0,$$
where
$$S_{\alpha}(t)=\frac{1}{2\pi i}\int_{\gamma} e^{\lambda t} \lambda^{\alpha-1}(\lambda^{\alpha}I-A)^{-1}d\lambda,
\quad t \ge t_0,$$ where the path $\gamma$ lies outside the sector
$\omega+S_{\theta}.$ If $S_{\alpha}(t)$ is integrable then the
solution is given by
$$x(t)=\int_{-\infty}^t S_{\alpha}(t-s)f(s,x(s),x_s)ds.$$
Now one can easily see that
$$z(t)-x(t)=S_{\alpha}(t)x_0-\int_t^{\infty}S_{\alpha}(s)f(t-s,x(t-s),x_{t-s}).$$
Hence for $f \in L^{p'}({\mathbb{R}^+\times X,X}), \ p' \in
[1,\infty)$ we have $v(t)-u(t) \rightarrow 0$ as $t\rightarrow
\infty.$
\begin{definition}
A function $x: [-r,T] \rightarrow X$ with $x(s)=\Phi(s)-g(x_t)(s), \quad t\in [0,T], s\in [-r,0]$ is said to be a mild solution to (\ref{meq}) if its restriction on $[0,T]$ satisfies
\begin{eqnarray}
x(t)&=& S_{\alpha}(t)(\phi(0)-g(x_t)(0)) + \int_{0}^tS_{\alpha}(t-s)Bu(s)ds\nonumber \\ &+&\int_{0}^tS_{\alpha}(t-s)f(s,x(s),x_s)ds,\quad \quad \forall t \in [0,t_1],\nonumber  \\
&=& S_{\alpha}(t)(\phi(0)-g(x_t)(0))+\sum_{0<t_k<t}S_{\alpha}(t-t_k)I_k(x(t_k^-)) \nonumber \\ &+& \int_{0}^tS_{\alpha}(t-s)Bu(s)ds+\int_{0}^tS_{\alpha}(t-s)f(s,x(s),x_s)ds,\quad \forall t \in (t_1,t_2],\nonumber  \\
&=& S_{\alpha}(t)(\phi(0)-g(x_t)(0))+\sum_{0<t_k<t}S_{\alpha}(t-t_k)I_k(x(t_k^-)) \nonumber \\ &+& \int_{0}^tS_{\alpha}(t-s)Bu(s)ds+\int_{0}^tS_{\alpha}(t-s)f(s,x(s),x_s)ds,\quad \forall t \in (t_2,t_3],\nonumber  \\
\vdots
\nonumber \\
&=& S_{\alpha}(t)(\phi(0)-g(x_t)(0))+\sum_{0<t_k<t}S_{\alpha}(t-t_k)I_k(x(t_k^-)) \nonumber \\ &+& \int_{0}^tS_{\alpha}(t-s)Bu(s)ds+\int_{0}^tS_{\alpha}(t-s)f(s,x(s),x_s)ds,\quad \forall t \in (t_K,T],\nonumber  \\
\end{eqnarray}
\end{definition}
\noindent Recently, Cuesta in \cite{custa}, theorem$1$, has proved that if
$A$ is a sectorial operator of type $\omega< 0$ for some $M > 0$ and $\theta \in [0,\pi(1-\frac{\alpha}{2}),$ then there exists $C> 0$ such that $$\|S_{\alpha}(t)\|_{B(X)} \le \frac{CM_1}{1+|\omega|t^{\alpha}}$$ for $t\ge0.$ For $t\in [0,T],$ it is easy to note that $$\|S_{\alpha}(t)\|_{B(X)} \le \frac{CM_1}{1+|\omega|t^{\alpha}} \le CM_1=M.$$

\begin{definition} (Definition 11.1 \cite{zeid})
Kuratowskii non-compactness measure: Let $D$ be a bounded set in metric space $(X,d)$, then Kuratowskii non-compactness measure, $\mu(D)$ is defined as $$\mu(D):=\inf \{\epsilon: \ \mbox{D covered by finitely many sets of diameters} \ \leq \epsilon\}.$$
\end{definition}

\begin{definition}(Definition 11.6 \cite{zeid})
Condensing map: Let $\Phi:D(\Phi) \subset X\rightarrow X$ be a bounded and continuous  operator on Banach space $X$ such that $\mu(\Phi(D))<\mu(D)$ for all bounded set $D\subset D(\Phi),$ provided $\mu(D)>0$, where $\mu$ is the Kuratowskii non-compactness measure, then $\Phi$ is called condensing map.
\end{definition}

\begin{definition} (Definition 2.9 \cite{zeid})
Compact map: A map $f:X\rightarrow X$ is said to be compact if $f$  is continuous and image of every bounded subset of $X$ under $f$ is pre-compact (closure is compact).
\end{definition}

\begin{theorem} (\cite{sado})
Let B be a convex, bounded and closed subset of a Banach space $X$ and $\Phi: B \rightarrow B$ be a condensing map. Then $\Phi$ has a fixed point in $B.$ \label{thm1}
\end{theorem}

\begin{lemma}(Example 11.7, \cite{zeid})
A map $ \Phi=\Phi_1+\Phi_2:X\rightarrow X$ is $k-$ contraction with $0\leq k< 1$ if
\begin{enumerate}
\item[(a)] $\Phi_1$ is $k-$ contraction i.e. $\|\Phi_1(x)-\Phi_1(y)\|_X \leq k\|x-y\|_X$ and
\item[(b)] $\Phi_2$ is compact,
\end{enumerate}
and hence $\Phi$ is a condensing map. \label{lem1}
\end{lemma}

\begin{definition}(Reachable set)
Let $x(T;u)$ be the state of (\ref{meq}) at $T$ with respect to admissible control $u$, then $\mathrm{R}(T)=\{x(T;u): u \in L^2([0,T],U)\}$ is called reachable set of (\ref{meq}).
\end{definition}

\begin{definition}(Approximate controllability)
System (\ref{meq}) is said to be approximate controllable if closure of the reachable set, $\overline{\mathrm{R}(T)}=X,$ i.e for any given $\epsilon >0$ the dynamic of the system can be steered from initial state $x_0$ to $\epsilon -$ neighbourhood of any of state in $X$ at time $T.$
\end{definition}

\begin{definition}(Solution space)
The solution space for our problem (\ref{meq}) is defined as follows
\begin{align*}
PC([-r,T],X) &= \mathcal{C}_r\cup PC(I,X), \ \mbox{where} \ \\
PC(I,X) &= \{ x:[0,T] \rightarrow X| x\in C([t_0,t_1],X)\cup C((t_k,t_{k+1}],X) k = 1,2 \cdots m,  \\ &  x(t_k^+) \ \mbox{and} \ x(t_k^-) \ \mbox{exist}, \ x(t_k)=x(t_k^-) \}
\end{align*}
is a Banach space with sup-norm $\|\cdot\|$, defined by $\|x\|=\sup \{\|x(t)\|_{X}: t\in I\}$ and $\mathcal{C}_r$ is an abstract phase space with norm $\|\cdot\|_r.$ Here solution space $PC([-r,T],X)$ also forms a Banach space sup-norm $\|\cdot\|$, defined by $\|x\|=\sup \{\|x(t)\|_{X}: t\in [-r,T]\}.$
\end{definition}

\section{\textbf{Approximate controllability}}
In this section, we first prove existence of mild solution of (\ref{meq}) in order to prove approximate controllability. For this we need the following assumptions:
\begin{enumerate}
\item[(A.1)] $f$ is bounded and Lipschitz, in particular $\|f(t,x,\phi)\|_X\leq M_1(t)$ and $\|f(t,x,\phi)-f(t,y,\psi)\|_X\leq L_1(t)\|x-y\|_X+L_2(t)\|\phi - \psi\|_r), \quad \forall \quad t \in [0,T], \quad x,y \in X, \quad \phi, \psi \in \mathcal{C}_r,$
\item[(A.2)] The function $f$ is continuous and compact in $I\times X\times \mathcal{C}_r$,
\item[(A.3)] $I_k$ is continuous, bounded and Lipschitz, in particular $\|I_k(x)\|_X\leq l_1$ and $\|I_k(x)-I_k(y)\|_X\leq l_2\|x-y\|_X, \quad k=1,2,\cdots ,K , \quad \forall x \in X,$
\item[(A.4)] $g$ is continuous and Lipschitz, in particular $\|g(\phi)-g(\psi)\|_X\leq l_g\|\phi-\psi\|_r,$
and $g(x_t)(s)=g(x(t+s)) \quad \forall t \in [0,T], \quad s \in [-r,0], \quad \phi , \psi \in \mathcal{C}_r ,$
\item[(A.5)] Define a linear operator $$\Gamma_0^T =\int_0^T S_{\alpha}(T-s)BB^*S_{\alpha}^*(T-s)ds$$ such that for each $h\in X, \quad z_{\epsilon}(h)=\epsilon(\epsilon I+\Gamma_0^TJ)^{-1}(h)$ converges to zero as $\epsilon \rightarrow 0$ in strong topology and $z_{\epsilon}(h)$ is a solution of $\epsilon z_{\epsilon}+\Gamma_0^T J(z_{\epsilon})=\epsilon h$, where $J:X\rightarrow X^*$ is defined as in \cite{mohmudov2},
\item[(A.6)] $$M(Kl_2+l_g+T^{1-p}(\|L_1\|_{\frac{1}{p}}+\|L_2\|_{\frac{1}{p}}))\Big(1
+\frac{MM_B^2T}{\epsilon}\Big)<1.$$ \label{thm2},
\end{enumerate}
where $L_1, L_2, M_1 \in L_{\frac{1}{p}}([0,T],\mathbb{R}^+),M_B=\|B\|_{B(X)}, 0<p<\alpha -1, l_1>0, l_2>0, l_g >0.$

\begin{remark}
Assumption $(A.5)$ is equivalent to approximate controllability of corresponding linear system of the system (\ref{meq}).
\end{remark}

\begin{remark}
If linear operator $B$ or solution operator $S_{\alpha}(t)$ is compact then controllability operator $Wu=\int_0^T S_{\alpha}(T-s)Bu(s)ds$ is also compact and it is not invertible if the Banach space $X$ is infinite dimensional. It has been established by Hernandez et al. \cite{hernandez} that exact controllability of semilinear evolution equations with the controllability operator $Wu=\int_0^T S_{\alpha}(T-s)Bu(s)ds$  valid if and only if $X$ is finite dimensional space.
\end{remark}

\begin{theorem}(\cite{mohmudov1})
Let $v:X\rightarrow X$ be a nonlinear operator and $x_{\epsilon}$ is solution of the following equation
\begin{eqnarray}
\epsilon z_{\epsilon} + \Gamma_0^T J(z_{\epsilon})=\epsilon v(\epsilon) \ \mbox{with}  \  \nonumber \\
\|v(z_{\epsilon})-w\|_X \ \mbox{as} \ \epsilon \rightarrow 0, \quad w \in X.
\end{eqnarray}
Then there exists a subsequence of $\{z_{\epsilon}\}$ strongly converges to zero as $\epsilon \rightarrow 0.$
\end{theorem} \label{thm3}

\begin{theorem}
Under the assumptions $(A.1)-(A.6)$ the fractional order semilinear evaluation equation (\ref{meq}) is approximately controllable.
\end{theorem} \label{thm4}

\textbf{Proof:}

With the help of assumption (A.5), for any arbitrary $x \in PC(I,X)$, define a control
\begin{eqnarray}
u_{\epsilon}(t,x)&=& B^*S_{\alpha}^*(T-t)J((\epsilon I+\Gamma_0^T J)^{-1}p(x)), \ \mbox{where} \nonumber \\
p(x)&=& h- S_{\alpha}(T)(\phi(0)-g(x_T)(0))-\sum_{0<t_k<T}S_{\alpha}(T-t_k)I_k(x(t_k^-)) \nonumber \\
&-& \int_{0}^TS_{\alpha}(T-s)f(s,x(s),x_s)ds. \nonumber
\end{eqnarray}
We break the proof in two parts in order to establish our result. \newline
Part-A: Existence of mild solution using the control $u_{\epsilon}(t,x)
[(t,x)$. \newline

Using the above control, we define the following operator
\begin{eqnarray}
Fx(t)&=& S_{\alpha}(t)(\phi(0)-g(x)(0))+\sum_{0<t_k<t}S_{\alpha}(t-t_k)I_k(x(t_k^-))\nonumber \\ &+& \int_{0}^tS_{\alpha}(t-s)Bu_{\epsilon}(s,x)ds+\int_{0}^tS_{\alpha}(t-s)f(s,x(s),x_s)ds,\quad \forall t \in [0,T]
\nonumber \\
\end{eqnarray}
and claim that it has a fixed point in $[0,T].$
\newline

Breaking the above function $F$ into two components, we get $$F_1x(t)=S_{\alpha}(t)(\phi(0)-g(x_t)(0))+\sum_{0<t_k<t}S_{\alpha}(t-t_k)I_k(x(t_k^-))+\int_{0}^tS_{\alpha}(t-s)Bu_{\epsilon}(s,x)ds$$ and $$F_2x(t)=\int_{0}^tS_{\alpha}(t-s)f(s,x(s),x_s)ds.$$

Now we prove the existence of mild solution using Sadovskii's fixed point theorem in four steps.
\newline
For any $q>0$ consider a closed ball $B_q=\{x\in PC(I,X)| \|x\|\leq q\}.$

\emph{\bf{Step-1:}} $F$ is self mapping on $B_q$ i.e. there exists $q>0$ such that $F(B_q)\subset B_q.$
\newline
Suppose that it is not true, then for each  $q>0$ there exists $x^q(\cdot)\in B_q $ such that $F(x^q) \notin B_q$
i.e. $F(x^q)(t)>q$.
\newline
Now
\begin{eqnarray}
q &<& \|F(x^q)(t)\|_X \nonumber \\
  & \leq & \|S_{\alpha}(t)\|_{B(X)}\|(\phi(0)-g(x_t^q)(0))\|_X+\sum_{0<t_k<t}\|S_{\alpha}(t-t_k)\|_{B(X)}\|I_k(x^q(t_k^-))\|_X \nonumber \\
  &+& \int_0^t \|S_{\alpha}(t-s)\|_{B(X)}\|Bu_{\epsilon}(s,x^q)\|_Xds +\int_0^t \|S_{\alpha}(t-s)\|_{B(X)}\|f(s,x^q(s),x_s^q)\|_Xds \nonumber \\
  &\leq & M\Big[|\phi(0)|+l_g\|x_t^q\|_r+\|g(0)\|_r+Kl_2\|x^q\| \nonumber \\ &+& \frac{M_B^2}{\epsilon}\int_0^t\Big(\|h-S_{\alpha}(T)(\phi(0)-g(x_T)(0))- \int_0^T S_{\alpha}(T-\tau)f(\tau,x^q(\tau,x_{\tau}^q)d\tau \|_X\Big)ds \nonumber \\ &+& \int_0^t \|f(s,x^q(s),x_s^q)\|_Xds\Big] \nonumber \\
  &\leq & M\Big[\Big(|\phi(0)|+|g(0)|+Kl_2+\int_0^T|f(s,0,0)|ds \nonumber \\ &+&\|x^q\|(Kl_2+l_g+T^{1-p}(\|L_1\|_{\frac{1}{p}}+\|L_2\|_{\frac{1}{p}}))\Big)(1+\frac{MM_B^2T}{\epsilon})\Big] . \nonumber
\end{eqnarray}
Dividing both sides by $q$ and taking limit as $q\rightarrow \infty ,$ we get
$$M\Big(Kl_2+l_g+T^{1-p}(\|L_1\|_{\frac{1}{p}}+\|L_2\|_{\frac{1}{p}}\Big)\Big(1+\frac{MM_B^2T}{\epsilon}\Big)\geq 1,$$
which is a contradiction and hence $F(B_q) \subset B_q$ i.e. $F$ is self mapping on $B_q.$
\newline

\emph{\bf{Step-2:}} $F_1$ is contraction, prove of this assertion is similar to the the above proof and hence we skip the overlapping steps. For $x, y \in B_q,$ we have
\begin{eqnarray}
&& \|F_1x(t)-F_1y(t)\|_X \nonumber \\ &&=\Big\|S_{\alpha}(t)(g(y)(0)-g(x)(0))+ \sum_{0<t_k<t}S_{\alpha}(t-t_k)I_k((x(t_k^-)-y(t_k^-)))
\nonumber \\
&&+ \int_0^t S_{\alpha}(t-s)B(u_{\epsilon}(s,x)-u_{\epsilon}(s,y))ds\Big\|_X \nonumber \\
&&\leq  M\Big[l_g+Kl_2+\frac{MM_B^2T}{\epsilon}\Big(l_g+Kl_2+T^{1-p}(\|L_1\|_{\frac{1}{p}}+\|L_2\|_{\frac{1}{p}})\Big)\Big]\|x-y\|.
\nonumber \\
\end{eqnarray}
Here $M\Big[l_g+Kl_2+\frac{MM_B^2T}{\epsilon}\Big(l_g+Kl_2+T^{1-p}(\|L_1\|_{\frac{1}{p}}+\|L_2\|_{\frac{1}{p}})\Big)\Big]<1,$
so $F_1$ is contraction.

\emph{\bf{Step-3}} $F_2$ is compact \newline
First we prove that $F_2$ is continuous. Let ${x^n}$ is sequence in $B_q$, then we  get the following
\begin{eqnarray}
&& f(s,x^n(s),x_s^n) \rightarrow f(s,x(s),x_s) \ \mbox{weakly and} \ \nonumber \\ && \|f(s,x^n(s),x_s^n)-f(s,x(s),x_s)\|_X\leq 2q(L_1(t)+L_2(t)) \nonumber
\end{eqnarray}
and hence using Lesbegue dominated convergence theorem, we see that $F_2$ is continuous. Now in the next step, we prove that $F_2$ is completely continuous.
\begin{eqnarray}
 && \|F_2x(\tau_2)-F_2x(\tau_1)\|_X \nonumber \\
&&\leq  MT^{2-\alpha}\Big[\int_0^{\tau_2}(\tau_2-s)^{\alpha-2}\|f(s,x(s),x_s)\|_Xds\nonumber \\ &&- \int_0^{\tau_1}(\tau_1-s)^{\alpha-2}\|f(s,x(s),x_s)\|_Xds\Big] \nonumber \\
&&\leq  MT^{2-\alpha}\Big[\int_0^{\tau_1}((\tau_1-s)^{\alpha-2}-(\tau_2-s)^{\alpha-2})\|f(s,x(s),x_s)\|_Xds \nonumber \\ &&+ \int_{\tau_1}^{\tau_2}(\tau_1-s)^{\alpha-2}\|f(s,x(s),x_s)\|_Xds\Big] \nonumber \\
&&\leq  MT^{2-\alpha}\Big[\int_0^{\tau_1}((\tau_1-s)^{\alpha-2}-(\tau_2-s)^{\alpha-2})M_1(s)ds \nonumber \\
&&+ \int_{\tau_1}^{\tau_2}(\tau_1-s)^{\alpha-2}M_1(s)ds\Big] \nonumber \\
&&\leq MT^{2-\alpha}\Big[\Big(\int_0^{\tau_1}((\tau_1-s)^{\alpha-2}-(\tau_2-s)^{\alpha-2})^{\frac{1}{1-p}}ds\Big)^{1-p}\Big(\int_0^{\tau_1}M_1^{\frac{1}{p}}(s)ds\Big)^p \nonumber \\ &&+\int_{\tau_1}^{\tau_2}((\tau_1-s)^{\alpha-2})^{\frac{1}{1-p}}ds\Big(\int_{\tau_1}^{\tau_2}M_1^{\frac{1}{p}}(s)ds\Big)^p\Big] \nonumber \\
&&\leq  MT^{2-\alpha}\Big[\Big(\int_0^{\tau_1}((\tau_1-s)^{\frac{\alpha-2}{1-p}}-(\tau_2-s)^{\frac{\alpha-2}{1-p}})ds\Big)^{1-p}\|M_1\|_{\frac{1}{p}} \nonumber \\ &&+ \Big(\int_{\tau_1}^{\tau_2}(\tau_2-s)^{\frac{\alpha-2}{1-p}}ds\Big)^{1-p}\|M_1\|_{\frac{1}{p}}\Big] \nonumber \\
&&\leq  MT^{2-\alpha}\|M_1\|_{\frac{1}{p}}(\tau_2-\tau_1)^{\alpha-p-1}.
\end{eqnarray}
The right side of the above expression is completely independent of $x$. Thus applying Arzela-Ascoli theorem for equicontinuous functions we conclude that $F_2(B_q)$ is relatively compact and hence $F_2$ is completely continuous on $I-\{t_1,t_2,t_3, \cdots , t_m\}.$ In similar way we can prove the equicontinuity of $F_2$ on $t=t_k^- \ \mbox{and} \ t=t_k^+ , k=1,2,3, \cdots , m.$ And thus $F_2$ is compact on $I=[0,T].$

\emph{\bf{Step-4:}} ($F$ is condensing) As $F=F_1+F_2, F_1$  is continuous and contraction, $F_2$ is compact, so lemma \ref{lem1} concludes that $F$ is condensing.
\newline
And hence by using Sadovskii's fixed point theorem \ref{thm1}, we conclude that (\ref{meq}) has a mild solution in $B_q.$

Part-B: Approximate controllability:
\newline
Let $x^{\epsilon}$ be a fixed point of $F$ in $B_q.$ Then we have
\begin{eqnarray}
u_{\epsilon}(t,x^{\epsilon}) &=& B^*S_{\alpha}^*(T-t)J((\epsilon I+\Gamma_0^T J)^{-1}p(x^{\epsilon})), where \nonumber \\
p(x^{\epsilon})&=& h-S_{\alpha}(T)(\phi(0)-g(x_T)(0))-\sum_{0<t_k<T}S_{\alpha}(T-t_k)I_k(x^{\epsilon}(t_k^-)) \nonumber \\ &-& \int_{0}^TS_{\alpha}(T-s)f(s,x^{\epsilon}(s),x_s^{\epsilon})ds, \nonumber
\end{eqnarray}
and satisfies the following equality
\begin{eqnarray}
x^{\epsilon}(T)& =& S_{\alpha}(T)(\phi(0)-g(x_T)(0))+\sum_{0<t_k<T}S_{\alpha}(T-t_k)I_k(x(t_k^-)) \nonumber \\
&+& \int_0^TS_{\alpha}(T-s)(f(s,x^{\epsilon}(s),x_s^{\epsilon})+Bu_{\epsilon}(s,x))ds \nonumber \\
&=& S_{\alpha}(T)(\phi(0)-g(x_T)(0))+\sum_{0<t_k<T}S_{\alpha}(T-t_k)I_k(x(t_k^-)) \nonumber \\
&+& \int_0^TS_{\alpha}(T-s)f(s,x^{\epsilon}(s),x_s^{\epsilon})ds+(\Gamma_0^TJ)((\epsilon I+\Gamma_0^T J)^{-1}p(x^{\epsilon})) \nonumber \\
&=& h-\epsilon(\epsilon I+\Gamma_0^T J)^{-1}p(x^{\epsilon})
\end{eqnarray}
From assumption $(A.1),$ we get
$$\int_0^T\|f(s,x^{\epsilon}(s),x_s^{\epsilon})\|_X^2ds \leq \int_0^TM_1^2(s)ds,$$
which follows that the sequence $\{f(\cdot,x^{\epsilon}(\cdot),x_{\cdot}^{\epsilon})\}$ is bounded. Using reflexivity of Banach space $X,$ there exists a subsequence, which for simplicity we denote by the same notation $\{f(\cdot,x^{\epsilon}(\cdot),x_{\cdot}^{\epsilon})\}$, converges weakly to a function $f(\cdot) \in L^2(I,X).$ Thus, we have
\begin{eqnarray}
\|v(x^{\epsilon})-w\|_X & =& \|\int_0^TS_{\alpha}(T-s)[f(s,x^{\epsilon}(s),x_s^{\epsilon})-f(s)]ds\|_X
\nonumber \\
 & = & \|\int_0^TS_{\alpha}(T-s)[f(s,x^{\epsilon}(s),x_s^{\epsilon})-f(s)]ds\|_X \nonumber \\ &&  \rightarrow 0 \quad \mbox{as} \ \epsilon \rightarrow 0,
\end{eqnarray}
by using the compactness of $f(\cdot),$ i.e., $f(\cdot) \rightarrow \int_0^{\cdot}S_{\alpha}(\cdot-s)f(s)ds,$ where $w=\int_0^TS_{\alpha}(T-s)f(s)ds.$ Hence by Theorem \ref{thm3}, we get $$\|x^{\epsilon}(T)-h\|_X=\|z_{\epsilon}\|_X \rightarrow 0$$ as $\epsilon \rightarrow 0$ for any $h \in X.$ Thus the approximate controllability the system (\ref{meq}) is established.

\noindent
\section{\textbf{Existence of optimal control }}
In this section we establish the existence of optimal control of the impulsive fractional
semilinear delay differential equations (\ref{meq}).
In order to proceed, let us define performance index $$J(u)=\int_0^TL(t,x(t),x_t,u(t))dt.$$
Our aim is to find a control $u^0 \in U_{ad}$ such that $J(u_0)\leq J(u)$ for all $u \in U_{ad},$
where $x$ denotes the mild solution of equation (\ref{meq}) corresponding to control $u \in U_{ad}.$ It is well known that the set of all admissible control $U_{ad}$ is a closed convex subset of $L_2(I,U).$
\newline
For existence of admissible control of problem (\ref{meq}), we assume the following:
\begin{enumerate}
\item[(L.1)] Functional $L$ is Borel measurable in $I\times X \times \mathcal{C}_r \times U$,
\item[(L.2)] $L(t,\cdot,\cdot,\cdot)$ is sequentially lower semi-continuous for almost all $t\in I$ in $X \times \mathcal{C}_r \times U,$
\item[(L.3)] $L(t,x,\phi,\cdot)$ is convex in $U,$
\item[(L.4)] $L(t,x,\phi,u)\geq (M_2(t)+d\|x\|_X+e\|\phi\|_r+f\|u\|_U^2), M_2\in L^1(I,\mathbb{R}_{+}).$
\end{enumerate}

\begin{theorem}
Under the assumptions of theorem (\ref{thm2}) and assumptions $(L.1)-(L.4),$ there exists an admissible optimal control of problem (\ref{meq}).
\end{theorem} \label{thm5}
\textbf{Proof:}
The main task is to minimize performance index $J(u).$ In order to prove that, we consider the following two cases:
\newline
\emph{Case-I} If $\inf \{J(u)| u\in U_{ad}\} =\infty $, then the result is obvious.
\newline
\emph{Case-II} If $\inf \{J(u)| u\in U_{ad}\} =\epsilon <\infty,$ then by definition of the infimum there exists a minimizing sequence ${u^n}\in U_{ad}$ such that $J(u^n)\rightarrow \inf \{J(u)| u\in U_{ad}\} \ \mbox{as} \ n\rightarrow \infty .$ As we know that $U_{ad}$ is a closed convex subset of reflexive Hilbert space $L_2(I,U)$, there exists a subsequence ${u^m}$ of the sequence ${u^n}$ which converges weakly to some point $u^0 \in U_{ad}.$
Corresponding to each $u^m$, there exists a mild solution $x^m$ of (\ref{meq})
i.e.
\begin{eqnarray}
x^m(t)&=& S_{\alpha}(t)(\phi(0)-g(x_t^m)(0)) + \int_{0}^tS_{\alpha}(t-s)Bu^m(s)ds\nonumber \\ &+&\int_{0}^tS_{\alpha}(t-s)f(s,x^m(s),x_s^m)ds,\quad \quad \forall t \in [0,t_1],\nonumber  \\
&=& S_{\alpha}(t)(\phi(0)-g(x_t^m)(0))+\sum_{0<t_k<t}S_{\alpha}(t-t_k)I_k(x^m(t_k^-)) \nonumber \\ &+& \int_{0}^tS_{\alpha}(t-s)Bu^m(s)ds+\int_{0}^tS_{\alpha}(t-s)f(s,x^m(s),x_s^m)ds,\quad \forall t \in (t_1,t_2],\nonumber  \\
&=& S_{\alpha}(t)(\phi(0)-g(x_t^m)(0))+\sum_{0<t_k<t}S_{\alpha}(t-t_k)I_k(x^m(t_k^-)) \nonumber \\ &+& \int_{0}^tS_{\alpha}(t-s)Bu^m(s)ds+\int_{0}^tS_{\alpha}(t-s)f(s,x^m(s),x_s^m)ds,\quad \forall t \in (t_2,t_3],\nonumber  \\
\vdots
\nonumber \\
&=& S_{\alpha}(t)(\phi(0)-g(x_t^m)(0))+\sum_{0<t_k<t}S_{\alpha}(t-t_k)I_k(x^m(t_k^-)) \nonumber \\ &+& \int_{0}^tS_{\alpha}(t-s)Bu^m(s)ds+\int_{0}^tS_{\alpha}(t-s)f(s,x^m(s),x_s^m)ds,\quad \forall t \in (t_K,T].\nonumber  \\
\end{eqnarray}
Similarly corresponding to $u^0,$ there exists a mild solution $x^0$ of (\ref{meq}) i.e.
\begin{eqnarray}
x^0(t)&=& S_{\alpha}(t)(\phi(0)-g(x_t^0)(0)) + \int_{0}^tS_{\alpha}(t-s)Bu^0(s)ds\nonumber \\ &+&\int_{0}^tS_{\alpha}(t-s)f(s,x^0(s),x_s^0)ds,\quad \quad \forall t \in [0,t_1],\nonumber  \\
&=& S_{\alpha}(t)(\phi(0)-g(x_t^0)(0))+\sum_{0<t_k<t}S_{\alpha}(t-t_k)I_k(x^0(t_k^-)) \nonumber \\ &+& \int_{0}^tS_{\alpha}(t-s)Bu^0(s)ds+\int_{0}^tS_{\alpha}(t-s)f(s,x^0(s),x_s^0)ds,\quad \forall t \in (t_1,t_2],\nonumber  \\
&=& S_{\alpha}(t)(\phi(0)-g(x_t^0)(0))+\sum_{0<t_k<t}S_{\alpha}(t-t_k)I_k(x^0(t_k^-)) \nonumber \\ &+& \int_{0}^tS_{\alpha}(t-s)Bu^0(s)ds+\int_{0}^tS_{\alpha}(t-s)f(s,x^0(s),x_s^0)ds,\quad \forall t \in (t_2,t_3],\nonumber  \\
\vdots
\nonumber \\
&=& S_{\alpha}(t)(\phi(0)-g(x_t^0)(0))+\sum_{0<t_k<t}S_{\alpha}(t-t_k)I_k(x^0(t_k^-)) \nonumber \\ &+& \int_{0}^tS_{\alpha}(t-s)Bu^0(s)ds+\int_{0}^tS_{\alpha}(t-s)f(s,x^0(s),x_s^0)ds,\quad \forall t \in (t_K,T].\nonumber  \\
\end{eqnarray}

We claim that $x^m$ converges strongly to $x^0.$

For $t\in [0,t_1]$, we have
\begin{eqnarray}
&& \|x^m(t)-x^0(t)\|_X \nonumber \\ && = \|S_{\alpha}(t)(g(x_t^0)(0)-g(x_t^m)(0))\|_X+\int_{0}^t\|S_{\alpha}(t-s)B(u^m(s)-u^0(s))\|_Xds \nonumber \\ && + \int_{0}^t\|S_{\alpha}(t-s)(f(s,x^m(s),x_s^m)-f(s,x^0(s),x_s^0))\|_Xds. \nonumber \\
&& \leq  M\Big[\|g(x_t^m)(0)-g(x_t^0)(0)\|_X + \int_{0}^t\|Bu^m(s)-Bu^0(s)\|_Xds \nonumber \\
&& + \int_{0}^t(L_1(s)\|x^m(s)-x(s)\|_X+L_2(s)\|x_s^m-x_s^0)\|_r)ds\Big].
 \label{eq1}
\end{eqnarray}
By taking the supremum of both side of equation (\ref{eq1}), we get
\begin{eqnarray}
\|x^m-x^0\| & \leq &  M\Big[\Big(l_g+T^{1-p}(\|L_1\|_{\frac{1}{p}}+\|L_2\|_{\frac{1}{p}})\Big)\|x^m-x^0\| \nonumber \\ &&+T\|Bu^m-Bu^0\|\Big].
\end{eqnarray}
For $t\in (t_1,t_2]$, we have
\begin{eqnarray}
&&\|x^m(t)-x^0(t)\|_X  \nonumber \\ &&= \|S_{\alpha}(t)(g(x_t^0)(0)-g(x_t^m)(0))\|_X+\sum_{0<t_k<t}\|S_{\alpha}(t-t_k)(I_k(x^m(t_k^-))\nonumber \\
 &&- I_k(x^0(t_k^-)))\|_X+\int_{0}^t\|S_{\alpha}(t-s)B(u^m(s)-u^0(s))\|_Xds \nonumber \\ &&+ \int_{0}^t\|S_{\alpha}(t-s)(f(s,x^m(s),x_s^m)-f(s,x^0(s),x_s^0))\|_Xds. \nonumber \\
&&\leq  M\Big[\|g(x_t^m)(0)-g(x_t^0)(0)\|_X+\|I_k(x^m(t_k^-))-I_k(x^0(t_k^-))\|_X\nonumber \\
 &&+ \int_{0}^t\|Bu^m(s)-Bu^0(s)\|_Xds + \int_{0}^t(L_1(s)\|x^m(s)-x(s)\|_X\nonumber \\
 &&+L_2(s)\|x_s^m-x_s^0)\|_r)ds\Big].
 \label{eq2}
\end{eqnarray}
By taking the supremum of both side of equation (\ref{eq2}), we get
\begin{eqnarray}
\|x^m-x^0\|  & \leq &  M\Big[\Big(l_g+l_2+T^{1-p}(\|L_1\|_{\frac{1}{p}}+\|L_2\|_{\frac{1}{p}})\Big)\|x^m-x^0\| \nonumber \\ &&+T\|Bu^m-Bu^0\|\Big].
\end{eqnarray}
For $t\in (t_K,T]$, we have
\begin{eqnarray}
&&\|x^m(t)-x^0(t)\|_X  \nonumber \\ &&= \|S_{\alpha}(t)(g(x_t^0)(0)-g(x_t^m)(0))+\sum_{0<t_k<t}\|S_{\alpha}(t-t_k)(I_k(x^m(t_k^-))\nonumber \\
 &&- I_k(x^0(t_k^-)))\|_X+\int_{0}^t\|S_{\alpha}(t-s)B(u^m(s)-u^0(s))\|_Xds \nonumber \\ &&+ \int_{0}^t\|S_{\alpha}(t-s)(f(s,x^m(s),x_s^m)-f(s,x^0(s),x_s^0))\|_Xds. \nonumber \\
&&\leq  M\Big[\|g(x_t^m)(0)-g(x_t^0)(0)\|_X+\|I_k(x^m(t_k^-))-I_k(x^0(t_k^-))\|_X\nonumber \\
 &&+ \int_{0}^t\|Bu^m(s)-Bu^0(s)\|_Xds + \int_{0}^t(L_1(s)\|x^m(s)-x(s)\|_X\nonumber \\
 &&+L_2(s)\|x_s^m-x_s^0)\|_r)ds\Big].
 \label{eq3}
\end{eqnarray}
Taking the supremum of both side of equation (\ref{eq3}), we get
\begin{eqnarray}
\|x^m-x^0\|  & \leq &  M\Big[\Big(l_g+Kl_2+T^{1-p}(\|L_1\|_{\frac{1}{p}}+\|L_2\|_{\frac{1}{p}})\Big)\|x^m-x^0\| \nonumber \\ &&+T\|Bu^m-Bu^0\|\Big].
\end{eqnarray}
Since $\Big(l_g+Kl_2+T^{1-p}(\|L_1\|_{\frac{1}{p}}+\|L_2\|_{\frac{1}{p}})\Big)<1$ and $\|Bu^m-Bu^0\|\rightarrow 0,$ it is easy to see that $x^m$ converges strongly to $x^0.$
\newline
Finally using Balder's theorem \cite{balder}, we get
\begin{eqnarray}
\epsilon &=& \lim_{m\rightarrow \infty}\int_0^TL(t,x^m(t),x_t^m,u^m(t))dt \nonumber \\
         &\leq & \int_0^TL(t,x^0(t),x_t^0,u^0(t))dt = J(u^0) \geq \epsilon .
\end{eqnarray}
Hence the result is followed that $J$ attains its minimum at $u^0 \in U_{ad}$.
\noindent
\section{\textbf{Example}} Consider the following fractional
relaxations – oscillation equation given by
\begin{eqnarray}
 && \frac{\partial^{\alpha}}{\partial t^{\alpha}} x(t, y) = \frac{\partial^2}{\partial y^2}x(t, y)
+\frac{\partial^{\alpha-1}}{\partial t^{\alpha-1}}(f_1(t,x(t,y))+ \int_{-r}^th(t-s)f_2(s,x(s,y))ds \nonumber \\ &&  +u(t,y)), \quad t \in I=[0,T],\ y \in \Omega=[0, \pi], \ \mbox{with boundary conditions} \ \nonumber \\
&& x(t,0)=x(t,\pi)=0 \quad t \in [0,T],\nonumber \\
&& \Delta x(t_k,y)=-x(t_k,y) \quad k=1,2, \cdots , K,\nonumber \\
&& x(s,y)+\sum_{i=1}^lc_ix(t_i+s,y)=\phi(s,y) \quad s\in [-r,0].
\label{eeq1}
\end{eqnarray}

Let $x(t)y=x(t,y)$ and assume $f(t,x(t),x_t)=f_1(t,x(t,y))+ \int_{-r}^th(t-s)f_(s,x(s,y)ds$ be a continuous function with respect to $t$ and satisfies Lipschitz condition in $x \ \mbox{and} \ x_t.$ Define the operator $Ax=
\frac{\partial^2x}{\partial y^2}-\omega x$ with domain $$D(A)=\{x \in L^2(0,\pi): x, x^{'} \ \mbox{are absolutely continuous and} \ x, x^{'}, x^{''} \in L^2(0,\pi)\}.$$ It is well known that for $\alpha =1,\ \mbox{sectorial operator,} \ A=\frac{\partial^2}{\partial y^2}-\omega$ generates an analytic semigroup and for $\alpha =2, \ \mbox{sectorial operator,} \ A=\frac{\partial^2}{\partial y^2}-\omega$ generates a cosine family of operators. Therefore the above problem can be posed as abstract problem (\ref{meq}) defined on $X=L^2(0,\pi)=U,$ and satisfies all the assumptions of Theorem \ref{thm4}. Hence the problem (\ref{eeq1}) has a mild solution in $[0,T]$ i.e. approximately controllable in $[0,T].$

Define performance index $$J(u)=\int_0^T\int_{\Omega}(\|x(t,y)\|_X^2+\|u(t,y)\|_X^2)dydt+\int_{-r}^0\int_{\Omega}\|x(t+s,y)\|_X^2dyds$$ with respect to problem (\ref{eeq1}).
Let $x(t)y=x(t,y)$ we see that
$$J(u)=\int_0^T(\|x(t)\|_X^2+\|x_t\|_r^2+\|u(t)\|_X^2)dt.$$
and satisfies all the assumptions of Theorem \ref{thm5}.  Hence there exists an admissible control $u^0 \in U_{ad} \in L^2(I,U)$
such that $J(u^0)\leq J(u), \quad \forall u  \in U_{ad}.$

%One can see that
%$$\int_0^{\infty}\frac{1}{1+|\omega|t^{\alpha}}=\frac{|\omega|^{\frac{-1}{\alpha}}\pi}{\alpha \sin
%\frac{\pi}{\alpha}}$$ for $\alpha \in (1,2).$ Thus the condition
%$2l\|S_{\alpha}\| <1$ reduces to
%$\frac{2L|\omega|^{\frac{-1}{\alpha}}\pi}{\alpha \sin
%\frac{\pi}{\alpha}}<1.$
%\textbf{monotone iterative technique}

\end{document}